\documentclass[12pt]{amsart}
\usepackage{amssymb}
\usepackage[colorlinks=true,pagebackref]{hyperref}
\usepackage[dvipsnames,usenames]{color}
\hypersetup{linkcolor=RawSienna,anchorcolor=BurntOrange,citecolor=OliveGreen,filecolor=BlueViolet,menucolor=Yellow,urlcolor=OliveGreen}

\renewcommand{\th}{{^\text{th}}}
\renewcommand{\phi}{\varphi}
\renewcommand{\epsilon}{\varepsilon}
\renewcommand{\theta}{\vartheta}
\newcommand{\mId}[1]{\mathcal{J}(#1)}

\newcommand{\rdwn}[1]{\lfloor #1 \rfloor}
\newcommand{\Oo}{\cO}

\def\ZZ{{\mathbb Z}}
\def\NN{{\mathbb N}}

\def\CC{{\mathbb C}}
\def\AA{{\mathbb A}}
\def\RR{{\mathbb R}}

\def\PP{{\mathbf P}}

\def\cJ{\mathcal{J}}

\def\cO{\mathcal{O}}

\def\fra{\mathfrak{a}}
\def\frb{\mathfrak{b}}

\def\frm{\mathfrak{m}}

\newcommand{\llbracket}{[\negthinspace[}
\newcommand{\rrbracket}{]\negthinspace]}

\DeclareMathOperator{\Hom}{Hom}

\newtheorem{lemma}{Lemma}[section]
\newtheorem{theorem}[lemma]{Theorem}
\newtheorem{corollary}[lemma]{Corollary}
\newtheorem{proposition}[lemma]{Proposition}

\theoremstyle{definition}
\newtheorem{definition}[lemma]{Definition}
\newtheorem{remark}[lemma]{Remark}

\newtheorem{example}[lemma]{Example}

\theoremstyle{remark}
\newtheorem*{remark*}{Remark}
\newtheorem*{note*}{Note}

\begin{document}

\title{Discreteness and rationality of $F$-thresholds}

\author[M.~Blickle]{Manuel~Blickle}
\address{FB Mathematik, Universit\"at  Duisburg-Essen, Standort Essen,
45117 Essen, Germany}
\email{{\tt manuel.blickle@uni-essen.de}}

\author[M. Musta\c{t}\v{a}]{Mircea~Musta\c{t}\v{a}}
\address{Department of Mathematics, University of Michigan,
Ann Arbor, MI 48109, USA}
\email{{\tt mmustata@umich.edu}}

\author[K.~E.~Smith]{Karen~E.~Smith}
\address{Department of Mathematics, University of Michigan,
Ann Arbor, MI 48109, USA}
\email{{\tt kesmith@umich.edu}}

\begin{abstract}
The $F$-thresholds are characteristic $p$ analogs of the jumping
coefficients for multiplier ideals in characteristic zero. We give
an alternative description of the $F$-thresholds of an ideal in a
regular and $F$-finite ring $R$, which enables us to settle two open
questions posed in \cite{MTW}. Namely,  we show that if the ring is,
in addition, essentially of finite type over a field, then the
$F$-thresholds are rational and discrete.
\end{abstract}

\thanks{The first author was partially supported by the \emph{DFG Schwerpunkt Komplexe Geometrie}.
The second and  third authors were partially supported by the NSF under grants DMS 0500127 and
 DMS 0500823, respectively.}

\maketitle

\markboth{M.~Blickle, M.~Musta\c t\v a and K.~E.~Smith}{Discreteness
and rationality of $F$-thresholds}

\section{Introduction}
In recent years,  multiplier ideals have played an important role in
higher dimensional birational geometry. For a given ideal $\fra$ on
a smooth variety $X$, and a real parameter $c>0$, the multiplier
ideal $\mId{\fra^c}$ is defined via a log resolution
of the pair $(X,\fra)$. Recall that a log resolution is a proper birational map $\pi: X' \rightarrow X$, with $X'$ smooth such that $\fra \Oo_{X'}$ defines a
simple normal crossing divisor $A=\sum_{i=1}^r a_iE_i$. Then, by definition,
\begin{equation}\label{multiplier_ideal}
    \mId{\fra^c} := \pi_*\Oo_{X'}(K_{X'/X}-\rdwn{cA}),
\end{equation}
and this is an ideal of $\Oo_X$ that does not depend on the chosen
log resolution. A jumping coefficient (also called a jumping number  or a jumping exponent)  of $\fra$ is a positive real
number $c$ such that $\mId{\fra^c} \neq \mId{\fra^{c-\epsilon}}$ for
every $\epsilon>0$. These invariants have been introduced and
studied in \cite{EinLazSmithVarolin.JumpingcoefficientsMultiplier}.
It follows from the formula~(\ref{multiplier_ideal}) that if $c$ is
a jumping coefficient, then $c\cdot a_i$ is an integer for some $i$.
In particular, every jumping coefficient is a rational number, and
the set of jumping coefficients of a given ideal is discrete.

Hara and Yoshida introduced in \cite{HY} a positive characteristic
analog of multiplier ideals, denoted by $\tau(\fra^c)$. This is a
generalized test ideal for a tight closure theory with respect to
the pair $(X,\fra^c)$. Similarly, one can define  jumping numbers
for these test ideals. These invariants were studied under the name
of $F$-thresholds in \cite{MTW}, where it was shown that they
satisfy many of the formal properties of the jumping coefficients in
characteristic zero.

We emphasize that the test ideals are not determined by a log
resolution of singularities, even in the cases where such a
resolution is known to exist. Instead, the definition uses the
Frobenius morphism and requires a priori infinitely many conditions
to be checked. This lack of built-in finiteness makes the question
of rationality and discreteness of the $F$-thresholds non-trivial,
and in fact these properties were left open in \cite{MTW}.

In the present paper we settle these questions in the case of a
regular ring $R$ that is essentially of finite type over an F-finite field:
more precisely, we show that for every ideal in such a ring, all
$F$-thresholds are rational and they form a discrete set. We mention
that in the case  when $\fra$ is a principal ideal in $R=k\llbracket
x,y\rrbracket$, where $k$ is a finite field, this statement has been
proved recently by Hara and Monsky in \cite{HM}. Our key result,
which implies discreteness and rationality of $F$-thresholds, is a
finiteness statement for test ideals in a polynomial ring. We show
in Proposition~\ref{prop6} below that if $\fra$ is an ideal in
$k[x_1,\ldots,x_n]$ that is generated in degree $\leq d$, then the
test ideal $\tau(\fra^c)$ is generated in degree $\leq \rdwn{c\cdot
d}$.

The structure of the paper is as follows. In Section
\ref{sec.BackGround} we recall the definition of test ideals and
$F$-thresholds in a slightly different setup from that in \cite{HY}
and \cite{MTW}. In fact,
 we prefer to work with a different definition of
test ideals that is more suitable for our purpose, and show later
that this is equivalent with the definition from \cite{HY}. Due to
this alternative definition, and to the slightly different setup, we
decided to develop the theory from scratch. We hope that this would
be beneficial for the reader, as with the present definition the
basic results become particularly transparent. After this initial
setup, we prove our main results in Section \ref{sec.DiscRat}.

\subsection{Acknowledgements}
We are grateful to Rob Lazarsfeld for useful discussions, to Craig
Huneke for a suggestion regarding the notation in the paper, and to
the referee for several helpful comments.

\section{Generalized test ideals and $F$-thresholds}\label{sec.BackGround}

In this section we collect some basic results about generalized test
ideals as introduced by Hara and Yoshida \cite{HY} (see also
\cite{Takagi.MultiplierTight} and \cite{HT}). We restrict ourselves
to working over an $F$-finite regular ring which allows us to use an
alternative description of these ideals. We subsequently prove that
our ideals agree with the ones defined in \cite{HY}. An advantage of
our definition is that the basic properties are very easy to prove.
We include some of the characteristic $p$ analogs of some of the
basic results on multiplier ideals, such as the Restriction Theorem
(Remark \ref{rem.restr}), Skoda's Theorem (Proposition \ref{prop7})
or the Subadditivity Theorem (Proposition \ref{lem4}(4)), all with very
elementary proofs.

We also recall the definition of $F$-thresholds from \cite{MTW},
though we do not restrict to the case of a \emph{local} ring as in
\emph{loc. cit}. In this more general framework we show that the set
of $F$-thresholds coincides with the set of jumping exponents for
the test ideals.

Let us fix the following notation: $R$ denotes a regular $F$-finite
ring of positive characteristic $p$. We stress here that $R$ is not
assumed to be local. The regularity of $R$ implies
 that the Frobenius morphism $F \colon R \to R$ sending
$x \in R$ to $x^p$ is flat \cite{Kunz}. This is equivalent to saying
that the functor $(F^e)^*$ of extending scalars via $F^e$ is exact.
By our assumption that $R$ is $F$-finite we just mean that $F$ is a
finite morphism. Since $F$ is also flat, it follows that $R$ is a
finitely generated locally free module over its subring $R^p$ of
$p\th$ powers.

For an ideal $J$ of $R$ and a positive integer $e$ we set
$J^{[p^e]}:=(u^{p^e}\vert u\in J)$. If $J=(u_1,\ldots,u_r)$, then
$J^{[p^e]}=(u_1^{p^e},\ldots,u_r^{p^e})$. Note that we have
$(J^{[p^e]})^{[p^{e'}]}=J^{[p^{e+e'}]}$. One easily verifies that
$(F^e)^*(R/J) \cong R/(J^{[p^e]})$. In particular, since $F$ is
(faithfully) flat we see that if $u$ is in $R$, then $u^p\in
J^{[p]}$ if and only if $u\in J$.

\begin{example}\label{rem0}
For a field $k$,  being $F$-finite precisely means that $[k\colon
k^p]<\infty$. One easily verifies, by explicitly giving a basis of
$R$ over $R^p$, that for such fields the polynomial ring
$k[x_1,\ldots,x_n]$ and even the power series ring
$k[[x_1,\ldots,x_n]]$ in finitely many variables are $F$-finite.
\newline
\indent More generally, if $A$ is an $F$-finite ring, then every
$A$-algebra that is essentially of finite type over $A$ is again
$F$-finite. To check this observe that if $S$ is a multiplicative
system in a ring $R$, and if $a_1,\ldots,a_N$ generate $R$ over
$R^p$, then $\frac{a_1}{1},\ldots,\frac{a_N}{1}$ generate $S^{-1}R$
over $(S^{-1}R)^p$.
\newline
\indent Conversely, if $k$ is a field and if $R$ is an $F$-finite
$k$-algebra essentially of finite type over $k$, then $[k\colon
k^p]<\infty$. Indeed, if $\frm$ is a maximal ideal of $R$ and
$K=R/\frm$, then $[K\colon K^p]<\infty$. Since $[K\colon k]<\infty$,
we deduce that $[k\colon k^p]<\infty$.
\newline
\indent If $(R,\frm)$ is a regular local ring that is $F$-finite,
then its completion $\widehat{R}$ is also $F$-finite (and regular).
In fact, since $R$ is local, $R$ is free over $R^p$. If
$a_1,\ldots,a_N$ give a basis of $R$ over $R^p$, we claim that these
elements also give a basis of $\widehat{R}$ over $\widehat{R}^p$.
Indeed, we have canonical isomorphisms
$$F^*(\widehat{R})=F^*(\projlim_{\ell}R/\frm^{\ell})
\simeq\projlim_{\ell}F^*(R/\frm^{\ell})\simeq\projlim_{\ell}
R/(\frm^{\ell})^{[p]}=\widehat{R}$$ (note that since $F$ is finite,
$F^*$ commutes with the above projective limit). Therefore the
Frobenius morphism on $\widehat{R}$ is obtained by base extension
from the Frobenius morphism on $R$, which implies our claim.
\end{example}

\subsection{The ideals $\frb^{[1/q]}$}
Let $\frb$ be an ideal of $R$ and let $q=p^e$. We now introduce certain
ideals $\frb^{[1/q]}$ which are crucial for our definition of the
generalized test ideals. More precisely, the test ideals of $\fra$
are constructed using all $(\fra^r)^{[1/p^e]}$. The ideals
$(\fra^r)^{[1/p^e]}$ were introduced in the case of a principal ideal
$\fra$ in \cite{AMBL} to study the $D$-module structure of a
localization of $R$. It was predicted in \cite{AMBL} that these
ideals might be related to tight closure theory and the present work
confirms this prediction. Furthermore, similar ideas also appear in
the description of generalized test ideals in \cite[Lemma 2.1]{HT}.

\begin{definition}
For an ideal $\frb$ of $R$ and $q=p^e$, where $e$ is a positive
integer, let $\frb^{[1/q]}$ denote the unique smallest ideal $J$ of
$R$ with respect to inclusion, such that
    \[
        \frb\subseteq J^{[q]}.
    \]
\end{definition}

 Since $R$ is a
finitely generated flat (hence projective) $R$-module via $F^e$, the
well-known lemma below implies that for every family of ideals $\{J_i\}_i$ of
$R$ we have $(\cap_iJ_i)^{[q]}=\cap_iJ_i^{[q]}$. Therefore
$\frb^{[1/q]}$ is well-defined. We make the convention that
$\frb^{[1/p^0]}=\frb$. Note in particular that $\frb^{[q'/q]}$ makes
sense for every $q=p^e$ and $q'=p^{e'}$.

\begin{lemma}\label{projective_module}
If $M$ is a finitely generated projective module over the Noetherian
ring $R$, then for every family $\{J_i\}_i$ of ideals of $R$ we have
$$\cap_i J_iM=\left(\cap_i J_i\right) M.$$
\end{lemma}

\begin{proof}
The assertion is clear when $M$ is free. The general case follows
since we can find a finitely generated free module $P$ such that $M$
is a direct summand of $P$.
\end{proof}

 The following lemma
collects some basic properties of the ideals $\frb^{[1/q]}$.

\begin{lemma}\label{lem1}
Let $\fra$ and $\frb$ be ideals in $R$ and $q=p^e$, $q'=p^{e'}$
where $e$ and $e'$ are positive integers.
\begin{enumerate}
\item[(1)] If $\fra\subseteq\frb$, then $\fra^{[1/q]}\subseteq
\frb^{[1/q]}$.
\smallskip
\item[(2)] $(\fra\cap\frb)^{[1/q]}\subseteq \fra^{[1/q]}
\cap \frb^{[1/q]}$ and $\fra^{[1/q]}+\frb^{[1/q]}=
(\fra+\frb)^{[1/q]}$.
\smallskip
\item[(3)] $(\fra\cdot\frb)^{[1/q]}\subseteq \fra^{[1/q]}\cdot
\frb^{[1/q]}$.
\smallskip
\item[(4)] $(\frb^{[q']})^{[1/q]}=\frb^{[q'/q]}\subseteq (\frb^{[1/q]})^{[q']}$.
\smallskip
\item[(5)] $\frb^{[1/qq']}\subseteq
(\frb^{[1/q]})^{[1/q']}$.
\smallskip
\item[(6)] $\frb^{[1/q]}\subseteq (\frb^{q'})^{[1/qq']}$.
\end{enumerate}
\end{lemma}

\begin{proof}
Statement (1) is straightforward from definition, and both  inclusions
"$\subseteq$" in (2) follow from (1). On the other hand, in order to
show that $(\fra+\frb)^{[1/q]}\subseteq\fra^{[1/q]}+\frb^{[1/q]}$,
it is enough to use the fact that
$$\fra+\frb\subseteq
\left(\fra^{[1/q]}\right)^{[q]}+\left(\frb^{[1/q]}\right)^{[q]}=
\left(\fra^{[1/q]}+\frb^{[1/q]}\right)^{[q]},$$ and the minimality
in the definition of $(\fra+\frb)^{[1/q]}$.

For (3), note that since we have $\fra\subseteq
(\fra^{[1/q]})^{[q]}$ and $\frb \subseteq (\frb^{[1/q]})^{[q]}$, we
deduce
$$\fra\cdot\frb\subseteq (\fra^{[1/q]})^{[q]}\cdot
(\frb^{[1/q]})^{[q]}=(\fra^{[1/q]}\cdot \frb^{[1/q]})^{[q]}.$$ The
inclusion in (3) then follows from the minimality in the definition
of $(\fra\cdot\frb)^{[1/q]}$.

Statement (4) is straightforward when $q' \geq q$, so we assume that $q = p^e \geq q' = p^{e'}$.
 For any ideal $J$ in a regular ring, we have $\frb^{[p^{e'}]}\subseteq J^{[p^{e}]}$ if and
only if $\frb\subseteq J^{[p^{e-e'}]}$. On the other hand, note that
since
$$\frb\subseteq
(\frb^{[1/p^e]})^{[p^e]}=\left((\frb^{[1/p^e]})^{[p^{e'}]}\right)^{[p^{e-e'}]},$$
we get $\frb^{[1/p^{e-e'}]}\subseteq (\frb^{[1/p^e]})^{[p^{e'}]}$ from the
minimality in the definition of $\frb^{[1/p^{e-e'}]}$.

Similarly, for (5) we use the minimality in the definition of
$\frb^{[1/qq']}$ and the fact that
$$\frb\subseteq (\frb^{[1/q]})^{[q]}\subseteq
\left((\frb^{[1/q]})^{[1/q']}\right)^{[qq']}.$$

Finally, for (6), note that (4) implies $\frb^{[1/q]}  =  (\frb^{[q']})^{[1/qq']}$, and this is contained in
$(\frb^{q'})^{[1/qq']}$ by (1).
\end{proof}

When $R$ is free over $R^{p^e}$ one has the following alternative
description of these newly defined ideals. Again, in the case of a
principal ideal this was already observed in \cite{AMBL}.

\begin{proposition}\label{prop1}
Suppose that $R$ is free over $R^q$ and let
 $e_1,\ldots,e_N$ be a basis of $R$ over $R^{q}$.
If $h_1,\ldots,h_s$ are  generators  of an ideal $\frb$ of $R$, and
if for every $i=1,\ldots, s$ we write
\begin{equation}\label{eq1}
h_i=\sum_{j=1}^Na_{i,j}^{q}e_j
\end{equation}
with $a_{i,j} \in R$, then $\frb^{[1/q]}=(a_{i,j} \mid i\leq s,j\leq
N)$.
\end{proposition}

Note that the proposition implies in particular that the description
therein does not depend on the chosen basis of $R$ over $R^{q}$ or
on the generators of $\frb$.

\begin{proof}[Proof of Proposition~\ref{prop1}]
It follows from (\ref{eq1}) that $(h_1,\ldots,h_s)\subseteq
(a_{i,j}^{q}\mid i\leq s,\, j\leq N)$ and therefore $\frb\subseteq
(a_{i,j}\vert i,j)^{[q]}$. Hence the inclusion ``$\subseteq$" in the
proposition follows immediately from the definition of
$\frb^{[1/q]}$.

For the reverse inclusion, suppose that $\frb\subseteq J^{[q]}$. If
$g_1,\ldots,g_m$ generate $J$, we may write
\begin{equation}\label{eq2}
h_i=\sum_{\ell=1}^mb_{\ell}g_{\ell}^{q}
\end{equation}
for some $b_{\ell}\in R$.  Remembering that $R$ is assumed free over $R^q$, consider the dual basis
$e_1^*,\ldots,e_N^*$ for $\Hom_{R^{q}}(R,R^{q})$, so
$e_i^*(e_j)=\delta_{i,j}$. It follows from (\ref{eq1}) that
$e_j^*(h_i)=a_{i,j}^{q}$. On the other hand, (\ref{eq2}) shows that
$$e_j^*(h_i)=\sum_{\ell}g_{\ell}^{q}e_j^*(b_{\ell})\in
J^{[q]}.$$ Therefore $a_{i,j}\in J$ for every $i$ and $j$, which
gives $(a_{i,j}\mid i\leq s,\, j\leq N)\subseteq \frb^{[1/q]}$.
\end{proof}

\begin{remark}
In \cite{AMBL} it is shown that the ideals $\frb^{[1/q]}$ have the
following $D$-module theoretic description. If $D^{(e)}_R$ denotes
the subring of the ring of differential operators $D_R$ which are
linear over $R^{q}$ (recall that $q=p^e$), then
$\left(\frb^{[1/q]}\right)^{[q]}$ is equal to the $D^{(e)}_R$-
submodule of $R$ generated by $\frb$. This interesting viewpoint
will not be exploited further here. It is however used successfully in our preprint \cite{BliMusSmi.HypF} to derive rationality and discreteness of $F$--thresholds of principal ideals in rings of formal power series, for example; thereby generalizing the 2 variable case treated in \cite{HM}.
\end{remark}

The following lemma shows that the formation of the ideals
$\frb^{[1/q]}$ commutes with localization and completion. In
particular, in order to compute $\frb^{[1/q]}$ we may always
localize so that $R$ is free over $R^q$ and then use
Proposition~\ref{prop1}.

\begin{lemma}\label{lem2}
Let $\frb$ be an ideal in $R$.
\begin{enumerate}
\item[(1)] If $S$ is a multiplicative system in $R$, then
$(S^{-1}\frb)^{[1/q]}=S^{-1}\left(\frb^{[1/q]}\right)$.
\item[(2)] If $R$ is local and $\widehat{R}$ is the completion of
$R$, then
$(\frb\widehat{R})^{[1/q]}=\left(\frb^{[1/q]}\right)\widehat{R}$.
\end{enumerate}
\end{lemma}

\begin{proof}
For (1), note that $\frb\subseteq \left(\frb^{[1/q]}\right)^{[q]}$
implies
$$S^{-1}\frb\subseteq
S^{-1}\left((\frb^{[1/q]})^{[q]}\right)=(S^{-1}\frb^{[1/q]})^{[q]}.$$
Therefore $(S^{-1}\frb)^{[1/q]}\subseteq
S^{-1}\left(\frb^{[1/q]}\right)$.

For the reverse inclusion, write $(S^{-1}\frb)^{[1/q]}=S^{-1}J$, for
some ideal $J$ such that $(J\colon s)=J$ for every $s\in S$. Using
the flatness of $F$ we see that $(J^{[q]}\colon s^{q})=J^{[q]}$ for
every $s\in S$, hence $(J^{[q]}\colon s)=J^{[q]}$. Since
$$S^{-1}\frb\subseteq (S^{-1}J)^{[q]}=S^{-1}(J^{[q]}),$$
it follows that $\frb\subseteq J^{[q]}$. Therefore
$\frb^{[1/q]}\subseteq J$, which gives ``$\supseteq$" in (1).

For (2) we will use Proposition~\ref{prop1}. Since $R$ is local, $R$
is free over $R^{q}$. Moreover, one shows as in Example~\ref{rem0}
that if $e_1,\ldots,e_N$ give a basis of $R$ over $R^{q}$, then
these elements also give a basis of $\widehat{R}$ over
$(\widehat{R})^{q}$, and  the assertion in (2) follows from
Proposition~\ref{prop1}.
\end{proof}

\subsection{Generalized test ideals}
We will now use the ideals $(\fra^r)^{[1/p^e]}$ to define the
generalized test ideals of Hara and Yoshida \cite{HY}. In
Proposition~\ref{prop2} it is shown that our definition coincides
with that in \cite{HY}.

\begin{lemma}\label{lem3}
Let $\fra$ be an ideal in $R$. If $r$, $r'$, $e$ and $e'$ are such
that $\frac{r}{p^e}\geq\frac{r'}{p^{e'}}$ and $e'\geq e$, then
\[
    (\fra^r)^{[1/p^e]}\subseteq (\fra^{r'})^{[1/p^{e'}]}.
\]
\end{lemma}
\begin{proof}
Note that $r'\leq rp^{e'-e}$, hence $\fra^{r'}\supseteq
\fra^{rp^{e'-e}}$. It follows from Lemma ~\ref{lem1} (1) and (6) that
$$(\fra^{r'})^{[1/p^{e'}]}\supseteq (\fra^{rp^{e'-e}})^{[1/p^{e'}]}
\supseteq (\fra^r)^{[1/p^e]}.$$
\end{proof}

Let $\fra$ be an ideal in $R$ and let $c$ be a positive real number.
If we denote by $\lceil x\rceil$ the smallest integer greater than or equal to $ x$, then
for every $e$ we have $\frac{\lceil cp^e\rceil}{p^e}
\geq\frac{\lceil cp^{e+1}\rceil}{p^{e+1}}$. It follows from
Lemma~\ref{lem3} that
\[
    \left(\fra^{\lceil cp^e\rceil}\right)^{[1/p^e]}\subseteq \left(\fra^{\lceil cp^{e+1}\rceil}\right)^{[1/p^{e+1}]}.
\]
\begin{definition}
With notation as above, the generalized test ideal of
$\fra$ with exponent $c$ is defined to be
\[
    \tau(\fra^c) = \bigcup_{e >0} \left(\fra^{\lceil cp^e\rceil}\right)^{[1/p^e]}.
\]
Since $R$ is Noetherian, this union stabilizes after
finitely many steps. In particular,  the test ideal
$\tau(\fra^c)$ is equal to  =$(\fra^{\lceil cp^e\rceil})^{[1/p^e]}$ for  all sufficiently large $e$.
\end{definition}

\begin{remark}\label{rem00}
We can write $R=R_1\times\cdots\times R_m$ where all $R_i$ are
($F$-finite) regular domains. An ideal $\fra$ in $R$ can be written
as $\fra=\fra_1\times\cdots\times\fra_m$ and it is clear that for
every $c$ we have
\begin{equation}\label{decomposition_product}
\tau(\fra^c)=\tau(\fra_1^c)\times\cdots\times\tau(\fra_m^c).
\end{equation}
 This allows us to assume that $R$ is a domain whenever this is
convenient.
\end{remark}

We make the convention that if $R$ is a domain, then
$\tau(\fra^0)=R$ unless $\fra=(0)$, in which case
$\tau(\fra^0)=(0)$. When $R$ is not necessarily a domain, we define
$\tau(\fra^0)$ such that the decomposition
(\ref{decomposition_product}) holds also when $c=0$.

\begin{proposition}\label{lem4}
Let $\fra$ and $\frb$ be  ideals of $R$.
\begin{enumerate}
\item[(1)] If $c_1<c_2$, then $\tau(\fra^{c_2})\subseteq
\tau(\fra^{c_1})$.
\item[(2)] If $\fra\subseteq\frb$, then $\tau(\fra^c)\subseteq
\tau(\frb^{c})$.
\item[(3)]
$\tau((\fra\cap\frb)^c)\subseteq\tau(\fra^c)\cap\tau(\frb^c)$ and
$\tau(\fra^c)+\tau(\frb^c)\subseteq\tau((\fra+\frb)^c)$.
\item[(4)]  $\tau((\fra\cdot\frb)^c)\subseteq
\tau(\fra^c)\cdot\tau(\frb^c)$.
\end{enumerate}
\end{proposition}

\begin{proof}
Lemma  ~\ref{lem1}(1) implies that $(\fra^{\lceil
c_2p^e\rceil})^{[1/p^e]} \subseteq (\fra^{\lceil
c_1p^e\rceil})^{[1/p^e]}$. By taking $e\gg 0$, we get the assertion
in (1). The other assertions follow similarly, by taking the limit
in the corresponding assertions from Lemma ~\ref{lem1}.
\end{proof}

\begin{remark}
The inclusion in (4) above is the analog  of the Subadditivity
Theorem for multiplier ideals in characteristic zero (see
\cite{positivity}, Theorem~9.5.20). See also Theorem~6.10 in
\cite{HY} for a different approach.
\end{remark}

A direct application of Lemma~\ref{lem2} (for $\frb=\fra^{\lceil
cp^e\rceil}$ with $e\gg 0$) shows that the formation of test ideals
commutes with localization and completion (see also \cite{HT}).

\begin{proposition}\label{prop4}
Let $\fra$ be an ideal in $R$ and $c$ a non-negative real number.
\begin{enumerate}
\item[(1)] If $S$ is a multiplicative system in $R$, then
$\tau((S^{-1}\fra)^c)=S^{-1}\tau(\fra^c)$.
\item[(2)] If $R$ is local and $\widehat{R}$ is the completion of
$R$, then $\tau((\fra\widehat{R})^c)=\tau(\fra^c)\widehat{R}$.
\end{enumerate}
\end{proposition}

We now  show that the family of test ideals $\tau(\fra^c)$ of a fixed
ideal $\fra$ is right continuous in $c$.

\begin{proposition}\label{prop3}
If $\fra$ is an ideal in $R$ and $c$ is a non-negative real number,
then there exists $\epsilon>0$ such that
$\tau(\fra^c)=(\fra^r)^{[1/p^e]}$ whenever
$c<\frac{r}{p^e}<c+\epsilon$. That is, $\tau(\fra^c) = \tau(\fra^{c'})$ where $c'$ is a rational number of the form $\frac{r}{p^e}$ which approximates $c$ from above  sufficiently well.
\end{proposition}

\begin{proof}
We show first that there is $\epsilon>0$ and an ideal $I$ in $R$
such that $(\fra^r)^{[1/p^e]}=I$ whenever
$c<\frac{r}{p^e}<c+\epsilon$. Indeed, otherwise we can find $r_m$
and $e_m$ for $m\geq 1$ such that $\frac{r_m}{p^{e_m}}$ form a
strictly decreasing sequence converging to $c$ and
$(\fra^{r_m})^{[1/p^{e_m}]}\neq (\fra^{r_{m+1}})^{[1/p^{e_{m+1}}]}$
for every $m$. After replacing this sequence by a subsequence we may
assume that $e_m\leq e_{m+1}$ for every $m$. By Lemma~\ref{lem3}, we
have $(\fra^{r_m})^{[1/p^{e_m}]}\varsubsetneq
(\fra^{r_{m+1}})^{[1/p^{e_{m+1}}]}$ for every $m$. Since this
sequence of ideals does not stabilize, we contradict the fact that
$R$ is Noetherian. Therefore we can find an ideal $I$ as claimed.

We show now that $I=\tau(\fra^c)$. By Remark \ref{rem00} we may
 assume that $R$ is a domain. If $\fra=(0)$, then
$(\fra^r)^{[1/p^e]}=(0)$ for every $r$ and $e$, and our assertion is
trivial. We assume henceforth that $\fra\neq (0)$. Let $e$ be large
enough such that $\tau(\fra^c)=\left(\fra^{\lceil
cp^e\rceil}\right)^{[1/p^e]}$ and $\frac{\lceil
cp^e\rceil}{p^e}<c+\epsilon$. If $cp^e$ is not an integer, then
$\frac{\lceil cp^e\rceil}{p^e}>c$, and we get $I=\tau(\fra^c)$.

Suppose now that $cp^e$ is an integer. After possibly replacing $e$
by a larger value, we may assume that $c+\frac{1}{p^e}<c+\epsilon$,
hence $I=\left(\fra^{cp^e+1}\right)^{[1/p^e]}\subseteq
\left(\fra^{cp^e}\right)^{[1/p^e]}=\tau(\fra^c)$. For the reverse
inclusion we need to show that $\fra^{cp^e} \subseteq I^{[p^e]}$.
Let $u \in \fra^{cp^e}$. If $e'\geq e$, then
$\frac{cp^{e'}+1}{p^{e'}}<c+\epsilon$, hence
$\fra^{cp^{e'}+1}\subseteq I^{[p^{e'}]}$. We deduce that if $v$ is a
nonzero element in $\fra$, then for every $e'\gg e$ we have
$vu^{p^{e'-e}}\in (I^{[p^e]})^{[p^{e'-e}]}$. This says that $u$ is
in the tight closure of the ideal $I^{[p^e]}$, which is equal to
$I^{[p^e]}$ since $R$ is a regular ring (see \cite{HH}). Therefore
we conclude that $u \in I^{[p^e]}$, and the proof is complete. [To avoid the explicit appearance
of tight closure one could first reduce to the local case and use
the freeness of $R$ over $R^{p^{e'-e}}$ to find, for $e' \gg e$, a
splitting $\phi: R \to R$ of the $(e'-e)\th$ iterate of the
Frobenius $F^{e'-e}: R \to R$ such that $\phi(v) = 1$. Applying
$\phi$ to an equation witnessing the membership $vu^{p^{e'-e}}\in
(I^{[p^e]})^{[p^{e'-e}]},$ it follows that $u \in
I^{[p^e]}$.]
\end{proof}

\begin{corollary}\label{cor0}
If $m$ is a positive integer, then for every $c \in \RR_{\geq 0}$
we have
$$\tau\left((\fra^m)^c\right)=\tau(\fra^{cm}).$$
\end{corollary}

\begin{proof}
It is clear that $\left((\fra^m)^r\right)^{[1/p^e]}=(\fra^{rm})^{[1/p^e]}$ for
every positive integer $r$ and $e$. Let $e$ be large enough such that
\[
    \tau\left((\fra^m)^c\right)=\left(\fra^{\lceil cp^e\rceil m}\right)^{[1/p^e]},\text{ and }\tau(\fra^{cm}) =\left(\fra^{\lceil cmp^e\rceil}\right)^{[1/p^e]}\, .
\]
If,  for some $e$ as above,  we have $cp^e\in\ZZ$, then our assertion is
clear. If this is not the case, then for $e\gg 0$,
 we have
$\frac{\lceil cp^e\rceil m}{p^e}$ larger than $cm$, but arbitrarily close to
$cm$ as $e$ gets large, and the statement follows from Proposition~\ref{prop3}.
\end{proof}

\begin{corollary}\label{cor1}
For every ideal $\fra$ in $R$ and every non-negative real number
$c$, there is $\epsilon>0$ such that $\tau(\fra^c)=\tau(\fra^{c'})$
for every $c'\in [c,c+\epsilon)$.
\end{corollary}

\begin{proof}
It is clear that we may take $\epsilon$ as given by
Proposition~\ref{prop3}.
\end{proof}

\begin{definition}
A positive real number $c$ is an \emph{F-jumping exponent} for
$\fra$ if $\tau(\fra^c)\neq \tau(\fra^{c-\epsilon})$ for every
positive $\epsilon$.
\end{definition}
Unless explicitly mentioned otherwise, we make the convention that
$0$ is also an $F$-jumping exponent. We will study the basic
properties of these numbers in the next section.

\begin{remark}\label{rem10}
Suppose that $K/k$ is an extension of perfect fields, and consider
the ring extension $R=k[x_1,\ldots,x_n]\subseteq
S=K[x_1,\ldots,x_n]$. If $\fra$ is an ideal in $R$, then
$\tau(\fra^c)\cdot S=\tau((\fra\cdot S)^c)$. Indeed, the monomials
of degree at most $p^e-1$ in each variable give a basis of both $R$
and $S$ over $R^{p^e}$ and $S^{p^e}$, respectively. It follows from
Proposition~\ref{prop1} that for every $r$ and $e$ we have
$(\fra^r)^{[1/p^e]}\cdot S=((\fra S)^r)^{[1/p^e]}$. In particular,
we see that $\fra$ and $\fra\cdot S$ have the same $F$-jumping
exponents.
\end{remark}

\subsection{The connection with the Hara-Yoshida definition}
We show now that the ideals $\tau(\fra^c)$ we have defined coincide with the
$\fra^c$-test ideals
introduced in \cite{HY},  temporarily denoted by $\tau'(\fra^c)$.
  To state Hara and Yoshida's definition, we need to recall their notion of
$\fra^c$-tight closure.

 For every $e \geq 1$, let $R^e$ denote the $R-R$-bimodule on
$R$, with the left $R$-module structure being the usual one, while
the right one being induced by $F^e\colon R\to R$. For any inclusion $N \subset M$ of $R$-modules,   there is an induced
map of left $R$-modules
\begin{equation*}
R^e \otimes N \rightarrow R^e\otimes M
\end{equation*}
whose image we denote by $N^{[p^e]}_M$.
An element $m \in M$ is said to be in the $\fra^c$-tight closure of $N$ in $M$ if there exists $u \in R$, not in any minimal prime, such that
\begin{equation*}
{\rm{image} } \,\, (u \fra^{\lceil cp^e\rceil} \otimes m) \in N^{[p^e]}_M
\end{equation*}
for all $e \gg 0$. The set of all such elements forms a submodule $N^{*\fra^c}_M$ of $M$.

 We are mainly interested in the case where $R$ is regular, $N$ is zero, and
 \begin{equation}
 M = E:=\bigoplus_{\frm \in {\rm {maxspec}} R}E_{R}(R/\frm)
 \end{equation}
 is the direct sum, over all the maximal ideals of $R$, of the
 injective hulls at the corresponding residue fields.With this notation,  Hara and Yoshida's definition can be stated as follows:

 \begin{definition} Let $R$ be a regular ring of prime characteristic $p$, $\fra$ an ideal of $R$ and $c$ any positive  real number.
 The $\fra^c$-test ideal, denoted $\tau'(\fra^c)$,  is the annihilator in $R$ of the $R$-module
 ${0_E^{*{\fra^c}}}$, where $E$ is  as above in (5).
 \end{definition}

 \begin{remark} More generally, Hara and Yoshida define the ideals
 $\tau'(\fra^c)$ even when $R$ is not regular as the annihilator of the
 $R$-module $\bigcup_M 0_M^{*\fra^c}$, where the union is taken over all finitely
 generated submodules of $E$.{\footnote{In the parlance of \cite{HH},
 this union could be called the {\it finitistic $\fra^c$-tight closure} of zero in $E$.}}
 In the regular case, however, it is easy to check that   $\bigcup_M 0_M^{*\fra^c} =  0^{*\fra^c}_E$,
 since the flatness of Frobenius implies $R^e \otimes M \subseteq
 R^e\otimes E$ for all submodules $M$ of $E$.
  \end{remark}

  \begin{remark} Hara and Yoshida's definition recovers the classical definition of the  test ideal,
  as introduced by Hochster and Huneke in \cite{HH},  by taking $\fra$ to be the unit ideal, and $c$ arbitrary; see
  \cite[Prop 8.23]{HH}. We caution the reader, however, that the name ``test ideal" is misleading if $\fra \neq (1)$, because the elements of $\tau'(\fra^c)$ are {\it not the test elements} for $\fra^c$-tight closure. See \cite{HY}, Definition 1.6 and Remark 1.10.
  \end{remark}

The next proposition says that these two notions of test ideals agree.
\begin{proposition}\label{prop2}
If $\fra$ is an ideal in a regular $F$-finite ring $R$, then for
every non-negative real number $c$ we have
$\tau(\fra^c)=\tau'(\fra^c)$.
\end{proposition}

 Before proving this, we recall a few well-known facts.
   First, if $(R, \frm)$ is a regular local ring of
 dimension $n$,
  then  the injective hull of the residue field  $E=E_R(R/\frm)$  can be identified with
the local
cohomology module $H_{\frm}^n(R)$. The Frobenius morphism on $R$
induces a Frobenius action on $H_{\frm}^n(R)$ and
hence also on $E$. This action, also denoted by $F$, can be
described as follows. If $x_1,\ldots,x_n$ form a minimal system of
generators of $\frm$ then
\[
    E\simeq R_{x_1\cdots x_n}/\sum_{i=1}^nR_{x_1\cdots \widehat{x_i}\cdots x_n}.
\]
Under this isomorphism, the Frobenius morphism on $E$ is induced by the Frobenius morphism on
$R_{x_1\cdots x_n}$.  In particular, an element  $\eta \in E$  represented by the fraction
$\frac{w}{(x_1\cdots x_n)^N}$ will be sent under Frobenius to the element $\eta^{p^e} \in E$ represented by the fraction
$$\frac{w^{p^e}}{(x_1\cdots x_n)^{Np^e}}.$$
Also using this representation of  $\eta \in E$, we have
${\rm {Ann}}_R \,\eta =   (x_1^N, \dots,  x_n^N):_R w.$
In particular, for all $\eta \in E$ we have
\begin{equation}\label{eq98}
{\rm{Ann}}_R \,(\eta^{p^e}) = ({\rm{Ann}}_R \,\eta)^{[p^e]},
\end{equation}
for all $e \geq 1$.
Indeed,
${\rm {Ann}}_R \, \eta^{p^e} =   (x_1^{p^eN}, \dots,  x_n^{p^eN}):_R w^{p^e} =  ((x_1^N, \dots,  x_n^{N}):_R w)^{[p^e]} ,$ with the last equality following from the flatness of Frobenius.

The Frobenius  action on $E$ induces an $R$-linear map $\nu^e\colon
R^e\otimes_RE\to E$ given by $\nu^e(a\otimes \eta)=a\eta^{p^e}$. In
the case that $R$ is regular, this map is  an isomorphism for all
$e\geq 1$. (See, for example, \cite{lyu}, Example~1.2 or \cite{bli},
p.20.) Thus we can conveniently compute  the $\fra^c$-tight closure
of zero in $E$ as
 the set of all $\eta \in E$ for which there exists $u \in R$ not in any minimal prime such that
 \begin{equation}\label{tc}
  u \fra^{\lceil c p^e \rceil}\eta^{p^e} = 0, \,\,\,\,\, {\rm{for \,\, all }}\,\,\, e \gg 0.
  \end{equation}

 We now prove two lemmas, both well-known to practitioners of tight closure, whose proofs we include to keep the paper self-contained. In the language of tight closure, the first says that, for a regular local ring, zero is tightly closure in the injective hull of the residue field.

\begin{lemma}\label{fact} With the above notation,
let $u$ be a non-zero element in a regular local ring $(R, \frm)$. If $\eta \in E = E_R(R/\frm)$ satisfies
  $u\eta^{p^e}=0$ for every $e\gg 0$, then $\eta=0$.
\end{lemma}

\begin{proof}
The hypothesis implies that
 $u\in {\rm Ann}_R(\eta^{p^e})$
 for all $e \gg 0$. But then by  (\ref{eq98}), we see that
 $u \in \left({\rm Ann}_R(\eta)\right)^{[p^e]}$
for every $e\gg 0$. By Nakayama's Lemma, since $u$ is nonzero, it follows that ${\rm
Ann}_R(\eta)=R$.
\end{proof}

The second lemma is essentially a special case of  the statement that $1$ is a $\fra^c$-test element in a regular local ring (see \cite[Theorem 1.7]{HY}).

\begin{lemma}\label{fact2}
If $R$ is a regular local ring, then
$$
0^{*\fra^c}_E = \big\{ \eta \in E \,\, | \,\, \fra^{\lceil cp^e \rceil} \eta^{p^e} = 0 \,\,\, {\rm{for\,\, all }} \,\, e \geq 1 \big\} .
$$
\end{lemma}

\begin{proof}
Clearly $\supseteq$ holds, so consider any $\eta \in 0^{*\fra^c}_E.$
We know that there is a non-zero element $u\in R$ such that
$u\fra^{\lceil cp^e\rceil}\eta^{p^e}=0$.

   Now, fix any $e' \geq 1$, and any $h \in \fra^{\lceil c p^{e'}\rceil}.$
    For all $e \geq 1$, we have
    $$
    h^{p^e} \in ( \fra^{\lceil c p^{e'}\rceil})^{p^e} \subseteq \fra^{\lceil c p^{e'+ e}\rceil}.
    $$
    Thus
    $$u ( h \eta^{p^{e'}})^{p^e} = u  h^{p^e} \eta^{p^{e'+e}} \in u  \fra^{\lceil c p^{e'+ e}\rceil} \eta^{p^{e + e'}},$$
which is zero for $e \gg 0$ by hypothesis.
By Lemma \ref{fact}, it follows that  $h \eta^{p^{e'}} = 0$ for all $e' \geq 1$.
Since $h$ was an arbitrary element of $\fra^{\lceil cp^{e'}\rceil}$, we conclude that   $\fra^{\lceil c p^{e'} \rceil} \eta^{p^{e'}} = 0$ for all $e' \geq 1$. The proof is complete.
\end{proof}

\begin{proof}[Proof of Proposition \ref{prop2} ]
The proof reduces immediately to the case $R$ is local. Indeed, since $E_R(R/\frm)$ is already a $R_{\frm}$-module, the $\fra^c$-tight closure of zero in $E_R(R/\frm)$ is the same as the   $(\fra R_{\frm})^c$-tight closure of zero in $E_{R_{\frm}}(R_{\frm}/\frm R_{\frm})$,
whence
\begin{equation}\label{eq99}
\tau'(\fra^c)=\cap_{\frm}\{f\in R\mid f/1 \in\tau\left((\fra
R_{\frm})^c\right)\}.
\end{equation}
This implies that  $\tau'(\fra^c)$ commutes with localization at maximal ideals, whereas Lemma \ref{lem2} implies
that $\tau(\fra^c)$ does.

Now, we assume $(R, \frm)$ is a regular local ring. To show that the ideals $\tau'(\fra^c)$  and
$\tau(\fra^c)$  coincide, it is enough, by Matlis duality, to show that their annihilators in $E$ coincide.
By definition (and Matlis duality), the annihilator of $\tau'(\fra^c)$ in $E$ is precisely $0^{*\fra^c}_E$.
It remains to  understand the annihilator in $E$ of $\tau(\fra^c)$.

By definition, $\tau(\fra^c) =   (\fra^{\lceil c p^e\rceil})^{[1/p^e]}$ for all $e \gg 0$. So an element
$\eta$ represented by a fraction
$\frac{w}{(x_1\cdots x_n)^N}$ is annihilated by  $\tau(\fra^c) $ if and only if
$$
(\fra^{\lceil c p^e\rceil})^{[1/p^e]}  \subseteq (x_1^N, \dots,
x_n^N): w.
$$
By definition of $\frb^{[1/q]}$, this holds if and only if
$$\fra^{\lceil c p^e\rceil}  \subseteq ((x_1^N, \dots, x_n^N): w)^{[p^e]} =  (x_1^{Np^e}, \dots, x_n^{Np^e}): w^{p^e},$$
 which in turn exactly says that
 $
\fra^{\lceil c p^e\rceil} \eta^{p^e} = 0$ for all $e \gg 0$.
We
conclude using Lemma \ref{fact2}  that the annihilator of $\tau(\fra^c)$ in $E$ is precisely $0^{*\fra^c}_E$. This completes the proof.
\end{proof}

\subsection{Skoda's Theorem}
For future reference, we include the following characteristic $p$
version of Skoda's Theorem, due to Hara and Takagi \cite{HT}. Since
the result in \emph{loc. cit.} is not stated in the generality we
will need, we include a proof for the benefit of the reader.

\begin{proposition}\label{prop7}
If $\fra$ is an ideal generated by $m$ elements, then for every
$c\geq m$ we have
$$\tau(\fra^c)=\fra\cdot\tau(\fra^{c-1}).$$
\end{proposition}

\begin{proof}
If $e$ is large enough, then $\tau(\fra^c)= \left(\fra^{\lceil
cp^e\rceil}\right)^{[1/p^e]}$ and
$\tau(\fra^{c-1})=\left(\fra^{\lceil
cp^e\rceil-p^e}\right)^{[1/p^e]}$. Therefore it is enough to show
that for every $r\geq mp^e$ we have
$$(\fra^r)^{[1/p^e]}=\fra\cdot (\fra^{r-p^e})^{[1/p^e]}.$$

The inclusion $\fra\cdot
\left(\fra^{r-p^e}\right)^{[1/p^e]}\subseteq (\fra^r)^{[1/p^e]}$
holds in fact for every $r\geq p^e$. Indeed, this says that
$\left(\fra^{r-p^e}\right)^{[1/p^e]}\subseteq
\left((\fra^r)^{[1/p^e]}\colon \fra\right)$, which is equivalent
with
\[
\fra^{r-p^e}\subseteq
\left((\fra^r)^{[1/p^e]}\colon\fra\right)^{[p^e]}
=\left(\left((\fra^r)^{[1/p^e]}\right)^{[p^e]}\colon
\fra^{[p^e]}\right).
\]
This holds since $\fra^{r-p^e}\cdot
\fra^{[p^e]}\subseteq\fra^r\subseteq
\left((\fra^r)^{[1/p^e]}\right)^{[p^e]}$.

Suppose now that $\fra=(h_1,\ldots,h_m)$. In order to prove the
reverse inclusion, note that if $r\geq m(p^e-1)+1$, then in the
product of $r$ of the $h_i$, at least one of these appears $p^e$
times. Therefore $\fra^r=\fra^{[p^e]}\cdot\fra^{r-p^e}$. We deduce
that
$$\fra^r\subseteq\fra^{[p^e]}\cdot\fra^{r-p^e}\subseteq\fra^{[p^e]}\cdot
\left(\left(\fra^{r-p^e}\right)^{[1/p^e]}\right)^{[p^e]}
=\left(\fra\cdot (\fra^{r-p^e})^{[1/p^e]}\right)^{[p^e]},$$ which
implies $(\fra^r)^{[1/p^e]}\subseteq\fra\cdot
\left(\fra^{r-p^e}\right)^{[1/p^e]}$.
\end{proof}

\begin{remark}\label{spread}
In Proposition \ref{prop7}, rather than taking $m$ to be the minimal number of generators for $\fra$, we can take $m$ to be the {\it analytic spread}  of $\fra$, that is, the minimal number of generators for any subideal with the same integral closure as $\fra. $ Indeed, this is immediate from the following lemma,
showing that the test ideals of $\fra$ depend only
on its  integral closure (see also \cite{HY} and
\cite{HT}).
\end{remark}

\begin{lemma}\label{closure}
If $\overline{\fra}$ denotes the integral closure of $\fra$, then
$\tau(\fra^c)=\tau(\overline{\fra}^c)$ for every $c$.
\end{lemma}

\begin{proof}
The inclusion $\tau(\fra^c)\subseteq\tau(\overline{\fra}^c)$ is
immediate from Proposition \ref{lem4} (2). For the reverse inclusion, note that by usual properties of
integral closure, there exists $m$ such that $\overline{\fra}^{m+\ell}
\subseteq\fra^{\ell}$ for every $\ell$. Corollary~\ref{cor1} gives
$c'>c$ such that $\tau(\fra^c)=\tau(\fra^{c'})$ and
$\tau(\overline{\fra}^c)=\tau(\overline{\fra}^{c'})$.

Using Corollary~\ref{cor0}, we see that
$$\tau(\overline{\fra}^{c'})=\tau((\overline{\fra}^{m+\ell})^{\frac{c'}{m+\ell}})
\subseteq\tau((\fra^{\ell})^{\frac{c'}{m+\ell}})=\tau(\fra^{c'-\frac{c'm}{m+\ell}}).$$
If $\ell\gg 0$, then $c<c'-\frac{c'm}{m+\ell}<c'$, hence by our
choice of $c'$ we get
$\tau(\overline{\fra}^c)\subseteq\tau(\fra^c)$.
\end{proof}

\begin{remark}\label{rem.restr}
If $\phi\colon R\to S$ is a morphism of regular, $F$-finite rings of
positive characteristic, and if $\frb$ is an ideal in $R$, then
$(\frb\cdot S)^{[1/p^e]}\subseteq \frb^{[1/p^e]}\cdot S$ for every
$e$. Indeed, since $\frb\subseteq
\left(\frb^{[1/p^e]}\right)^{[p^e]}$ we have $\frb\cdot S\subseteq
\left(\frb^{[1/p^e]}\cdot S\right)^{[p^e]}$.
\newline \indent
We deduce that for every non-negative $c$ we have
 $\tau((\fra\cdot S)^c)\subseteq
\tau(\fra^c)\cdot S$. This is an analogue of the Restriction Theorem
for multiplier ideals in characteristic zero (see \cite{positivity},
Examples~9.5.4 and 9.5.8). For a different argument in a more
general (characteristic $p$) framework, see \cite{HY}, Theorems~4.1
and 6.10.
\end{remark}

\subsection{$F$-jumping exponents and $F$-thresholds} In \cite{MTW} $F$-jumping exponents
of an ideal $\fra$ were described as $F$-thresholds. Since the statements and the proofs in \emph{loc. cit.} were given in
the local case, we review them here for the reader's convenience.

Let $\fra$ be an ideal in $R$. For a fixed ideal $J$ in $R$ such
that $\fra\subseteq {\rm rad}(J)$ and for an integer $e > 0$ we
define $\nu^J_{\fra}(p^e)$ to be the largest non-negative integer
$r$ such that $\fra^r\not\subseteq J^{[p^e]}$ (if there is no such
$r$, then we put $\nu^J_{\fra}(p^e)=0$). If $\fra^r\not\subseteq
J^{[p^e]}$, then $\fra^{pr}\not\subseteq J^{[p^{e+1}]}$. Indeed,
otherwise we get $(\fra^r)^{[p]}\subseteq J^{[p^{e+1}]}$, hence
$\fra^r\subseteq J^{[p^e]}$, a contradiction. Therefore
$$\frac{\nu^J_{\fra}(p^e)}{p^e}\leq
\frac{\nu^J_{\fra}(p^{e+1})}{p^{e+1}},$$ hence we may define the
\emph{$F$-threshold} of $\fra$ with respect to $J$ as
$$c^J(\fra):=\lim_{e\to\infty}\frac{\nu^J_{\fra}(p^e)}{p^e}
=\sup_{e\geq 1}\frac{\nu^J_{\fra}(p^e)}{p^e}.$$

Note that if $\fra$ is generated by $s$ elements, then
$\fra^{s(p^e-1)+1}\subseteq\fra^{[p^e]}$. If $\fra^{\ell}\subseteq
J$, then $\nu^J_{\fra}(p^e)\leq \ell(s(p^e-1)+1)-1$ for every $e$.
Therefore $c^J(\fra)\leq s\ell$, in particular $c^J(\fra)$ is
finite.

The following proposition relates the $F$-thresholds with the
generalized test ideals of $\fra$.

\begin{proposition}\label{prop5}
Let $\fra$ be an ideal in $R$.
\begin{enumerate}
\item[(1)] If $J$ is an ideal in $R$ such that $\fra\subseteq
{\rm rad}(J)$, then
$$\tau(\fra^{c^J(\fra)})\subseteq J.$$
\item[(2)] If $c$ is a non-negative real number, then
$\fra\subseteq {\rm rad}(\tau(\fra^c))$ and
$$c^{\tau(\fra^c)}(\fra)\leq c.$$
\end{enumerate}
\end{proposition}

\begin{proof}
For (1), note that by Corollary~\ref{cor1} there is $c'>c^J(\fra)$
such that $I:=\tau(\fra^{c^J(\fra)})=\tau(\fra^{c'})$. Suppose now
that $e\gg 0$, so $\tau(\fra^{c'})=\left(\fra^{\lceil
c'p^e\rceil}\right)^{[1/p^e]}$.

Since $c'>c^J(\fra)$ and $e$ is large enough, we have $\lceil
c'p^e\rceil\geq \nu^J_{\fra}(p^e)+1$, hence
$$\fra^{\lceil c'p^e\rceil}\subseteq J^{[p^e]}.$$
This implies that $I\subseteq J$, as required.

For (2), let $e$ be large enough, so
$\tau(\fra^c)=\left(\fra^{\lceil cp^e\rceil}\right)^{[1/p^e]}$. By
definition, we have $\fra^{\lceil cp^e\rceil}\subseteq
\tau(\fra^c)^{[p^e]}$, which implies
$$\nu^{\tau(\fra^c)}_{\fra}(p^e)\leq \lceil cp^e\rceil-1.$$
Dividing by $p^e$ and letting $e$ go to infinity, we get the
required inequality.
\end{proof}

\begin{corollary}\label{cor2}
For every ideal $\fra$ in $R$, the set of $F$-jumping exponents for
$\fra$ is equal to the set of $F$-thresholds of $\fra$ (as we range over
all possible ideals $J$).
\end{corollary}

\begin{proof}
We show first that if $\alpha$ is an $F$-jumping exponent for
$\fra$, then $\alpha=c^J(\fra)$ for $J=\tau(\fra^{\alpha})$. Indeed,
Proposition~\ref{prop5}(2) gives $c^J(\fra)\leq\alpha$. Therefore,
by Proposition \ref{lem4}(1),
$J=\tau(\fra^{\alpha})\subseteq\tau(\fra^{c^J(\fra)})$. We also have
the reverse inclusion by Proposition~\ref{prop5}(1). Since $\alpha$
is an $F$-jumping exponent, we must have $\alpha=c^J(\fra)$.

Suppose now that $\alpha=c^J(\fra)$ for some $J$ whose radical contains $\fra$.  We need to show that
$\alpha$ is a jumping exponent. If this is not the case, then there
is $\alpha'<\alpha$ such that
$\tau(\fra^{\alpha})=\tau(\fra^{\alpha'})$. Using
Proposition~\ref{prop5}(1) we get $\tau(\fra^{\alpha'}) \subseteq
J$. If $e$ is large enough, then
$\tau(\fra^{\alpha'})=\left(\fra^{\lceil
\alpha'p^e\rceil}\right)^{[1/p^e]}$. Therefore
$\fra^{\lceil\alpha'p^e\rceil} \subseteq J^{[p^e]}$, hence
$\nu^J(\fra)\leq\lceil \alpha'p^e\rceil -1$. Dividing by $p^e$ and
letting $e$ go to infinity, we get $c^J(\fra)\le\alpha'$, a
contradiction. This completes the proof of the corollary.
\end{proof}

\subsection{Mixed generalized test ideals} Using our approach it is easy to
define also "mixed" test ideals, as in \cite{HY} and \cite{HT}. For
example, suppose that $\fra_1,\ldots,\fra_r$ are ideals in the
regular $F$-finite ring $R$, and let $c_1,\ldots,c_r\in\RR_+$. For
every $e\geq 1$, consider
$$I_e:=\left(\fra_1^{\lceil
c_1p^e\rceil}\cdot\ldots\cdot\fra_r^{\lceil c_r
p^e\rceil}\right)^{[1/p^e]}.$$ As before, one can check easily that
$I_e\subseteq I_{e+1}$ for every $e$, and since $R$ is Noetherian,
this sequence of ideals stabilizes. Its limit is the test ideal
$\tau(\fra_1^{c_1}\cdots\fra_r^{c_r})$. It is straightforward to
check that the basic properties we have discussed so far generalize
to this setting. For example, the Subadditivity Formula from
Proposition ~\ref{lem4} can be generalized with the same proof to the
following statement: if $\fra_1,\ldots,\fra_r$ are ideals in $R$,
and if $\lambda_1,\ldots,\lambda_r\in\RR_+$, then
$$\tau(\fra_1^{\lambda_1}\cdots\fra_r^{\lambda_r})\subseteq\tau(\fra_1^{\lambda_1})\cdot
\ldots\cdot\tau(\fra_r^{\lambda_r})$$ (see also \cite{HY},
Theorem~6.10). Similarly, the Summation Formula due to Takagi (see
\cite{Tak1}, Theorem~3.1) can be easily proved using our definition.

\section{Discreteness and rationality}\label{sec.DiscRat}
In this section we prove our main result.
\begin{theorem}\label{thm1}
Let $k$ be a field of characteristic $p>0$ and let $R$ be a regular
$F$-finite ring, essentially of finite type over $k$.  Suppose that
$\fra$ is an ideal in $R$.
\begin{enumerate}
\item[(1)] The set of $F$-jumping exponents of
$\fra$ is discrete (in every finite interval there are only finitely
many such numbers).
\item[(2)] Every $F$-jumping exponent of $\fra$ is a rational
number.
\end{enumerate}
\end{theorem}

We will reduce the proof of the theorem to the case
$R=k[x_1,\ldots,x_n]$. We start with some preliminary results. The
first proposition, of independent interest, gives an effective bound
for the degrees of the generators of the ideals $\tau(\fra^c)$ in
terms of the degrees of the generators of $\fra$. It is our main
ingredient for the proof of the theorem in the polynomial ring case.
For a real number $t$ we will denote by $\lfloor t\rfloor$ the
largest integer less than or equal to $ t$, and by $\{t\}$ the fractional part
$t-\lfloor t\rfloor$.

\begin{proposition}\label{prop6}
Let $\fra$ be an ideal in the polynomial ring $k[x_1,\ldots,x_n]$,
where $k$ is a field of characteristic $p$ such that $[k\colon
k^p]<\infty$. If $\fra$ can be generated by polynomials of degree at
most $d$, then for every non-negative real number $c$, the ideal
$\tau(\fra^c)$ can be generated by polynomials of degree at most
$\lfloor cd\rfloor$.
\end{proposition}

\begin{proof}
Fix first $r$ and $e$. The ideal $\fra^r$ is generated by
polynomials of degree at most $rd$. Choose such generators
$h_1,\ldots,h_s$ for $\fra^r$.

Let $b_1,\ldots,b_m$ be a basis of $k$ over $k^{p^e}$, and consider
the basis of $R=k[x_1,\ldots,x_n]$ over $R^{p^e}$ given by
$$\{b_ix^u\mid i\leq m, u\in\NN^n,0\leq u_j\leq p^e-1\,{\rm for}\,{\rm all}\,j\}$$
(if $u=(u_1,\ldots,u_n)$, then we put $x^u=x_1^{u_1}\cdots
x_n^{u_n}$). If we write
$$h_j=\sum_{i,u}a_{i,u}^{p^e}b_ix^u,$$
with $a_{i,u}\in R$, then for every $i$ and $u$ we have
$\deg(a_{i,u}^{p^e})\leq \deg(h_j)\leq rd$. It follows from
Proposition~\ref{prop1} that $(\fra^r)^{[1/p^e]}$ can be generated
by polynomials of degree at most $\frac{rd}{p^e}$. But by
 Proposition \ref{prop3}, $\tau(\fra^c) = (\fra^r)^{[1/p^e]}$ for some $\frac{r}{p^e} \geq c$ closely approximating $c$.
 Thus for sufficiently close approximations $\frac{r}{p^e}$,
 $\tau(\fra^c) $ is generated by polynomials of degree at most $\lfloor \frac{rd}{p^e}\rfloor =
 \lfloor cd \rfloor.$ The proof is complete.
\end{proof}

\begin{remark}\label{char0}
As Rob Lazarsfeld points out, the analog of Proposition \ref{prop6}
 for multiplier ideals can be deduced from a vanishing statement. Indeed,
suppose that $\fra\subseteq\CC[x_1,\ldots,x_n]$ is an ideal
generated in degree at most $d$, and denote by $\cJ(\fra^c)$ its
multiplier ideal. After homogenizing $\fra$ and taking the
associated sheaf on $\PP^n$, we get a sheaf of ideals
$\widetilde{\fra}$ whose restriction to $\AA^n=(x_0\neq
0)\subseteq\PP^n$ is $\fra$, and such that $\widetilde{\fra}\otimes
\cO_{\PP}(d)$ is globally generated. If $H$ denotes the hyperplane
class on $\PP^n$, then $\lfloor cd\rfloor H\sim
K_{\PP^n}+(n+1+\lfloor cd\rfloor)H$. For $1\leq i\leq n$ we have
$(n+1+\lfloor cd\rfloor-i-cd)H$ ample, and therefore
$$H^i(\PP^n,\cJ(\widetilde{\fra}^c)\otimes\cO_{\PP}(\lfloor
cd\rfloor-i))=0$$ by the Nadel Vanishing Theorem (see
\cite{positivity}, Corollary~9.4.15). Therefore
$\cJ(\widetilde{\fra}^c)\otimes\cO_{\PP}(\lfloor cd\rfloor)$ is
$0$-regular in the sense of Castelnuovo-Mumford regularity, hence it
is generated by global sections. This implies that $\cJ(\fra^c)$ is
generated in degree at most $\lfloor cd\rfloor$.
\end{remark}

\begin{proposition}\label{lem5}
Let $\fra$ be an ideal in a regular, $F$-finite ring $R$.
\begin{enumerate}
\item[(1)] If $\alpha$ is an $F$-jumping exponent for $\fra$, then
$p\alpha$ is an $F$-jumping exponent, too.
\item[(2)] If $\fra$ can be generated by $m$ elements, and if
$\alpha>m$ is an $F$-jumping exponent for $\fra$, then $\alpha-1$ is
an $F$-jumping exponent, too.
\end{enumerate}
\end{proposition}

\begin{remark}
For a stronger statement in (2), we can take $m$ to be any integer at least as large as  the analytic spread of $\fra$. (This follows from Remark \ref{spread}.) For example, if $R$ is local with an infinite residue field, we can take $m = \dim R$.
\end{remark}

\begin{proof}
For (1), note that by Corollary~\ref{cor2} there is an ideal $J$
containing $\fra$ in its radical such that $\alpha=c^J(\fra)$. It is
clear that $\nu^{J^{[p]}}_{\fra}(p^e)=\nu^J_{\fra}(p^{e+1})$, so
$c^{J^{[p]}}(\fra)=p\cdot c^J(\fra)$. Hence $p\alpha$ is an
$F$-jumping exponent of $\fra$ by Corollary~\ref{cor2}.

For (2), suppose that $\alpha-1$ is not an $F$-jumping exponent. Let
$\epsilon>0$ be such that
$\tau(\fra^{\alpha-1})=\tau(\fra^{\alpha-1-\epsilon})$ and
$\alpha-\epsilon>m$. It follows from Proposition~\ref{prop7} that
$\tau(\fra^{\alpha})=\tau(\fra^{\alpha-\epsilon})$, hence $\alpha$
is not an $F$-jumping exponent, a contradiction.
\end{proof}

The following proposition relates the generalized test ideals of two
different ideals defining the same scheme. For the characteristic
zero analogue in the context of multiplier ideals, see
Proposition~2.3 in \cite{Mu}.

\begin{proposition}\label{lem6}
Let $R$ be a regular, $F$-finite ring of positive characteristic,
and $I$ an ideal in $R$ of pure codimension $r$, such that $S=R/I$
is regular. If $\fra$ is an ideal in $R$ containing $I$, then for
every non-negative real number $c$ we have
$$\tau((\fra/I)^c)=\tau(\fra^{c+r})\cdot S.$$
\end{proposition}

\begin{proof}
By Proposition~\ref{prop4}, forming generalized test ideals commutes
with localization and completion. Therefore it is enough to prove
the case where  $R$ is local and complete. Since $I$ is generated by
part of a regular system of parameters for $R$, by doing induction
on $r$ we see that it is enough to prove the case $r=1$. Therefore
we may assume that $R=k\llbracket x_1,\ldots,x_n\rrbracket$ and
$I=(x_n)$. Note that $[k\colon k^p]<\infty$, since $R$ is $F$-finite.

We claim that
\begin{equation}\label{claim}
\left(\fra^{r+p^e}\right)^{[1/p^e]}\cdot
S=\left((\fra/(x_n))^{r +1}\right)^{[1/p^e]}
\end{equation} for any integer $r$ and all $e \gg 0$. This implies the desired statement.
 Indeed, if $e$ is large enough, then we have
$$\tau(\fra^{c+1})\cdot S=\left(\fra^{\lceil cp^e\rceil+p^e}\right)^{[1/p^e]}\cdot S=
\left((\fra/(x_n))^{\lceil
cp^e\rceil+1}\right)^{[1/p^e]}=\tau((\fra/(x_n))^c).$$ The last
equality follows from Proposition ~\ref{prop3},  since $\frac{\lceil cp^e\rceil +1}{p^e}$ is larger
than $c$, but  arbitrarily close to $c$ as $e$ gets very large.

We prove now the claim (\ref{claim}), using the description in Proposition~\ref{prop1}. Let
$a_1,\ldots,a_m$ be a basis of $k$ over $k^{p^e}$, and consider the basis of $R$
over $R^{p^e}$ given by
\[
\{a_ix^u\mid i\leq m,\, u=(u_j)\in\NN^n, 0\leq u_j\leq p^e-1\,{\rm
for}\,{\rm all}\,j\}.
\]
Write $\fra=(x_n)+\frb$, where $\frb$ is generated by power series in
$k\llbracket x_1,\ldots,x_{n-1}\rrbracket$. We have
$$\fra^{r+p^e}=\sum_{i=0}^{r+p^e}x_n^i\frb^{r+p^e-i}.$$
The generators of $\left(\fra^{r+p^e}\right)^{[1/p^e]}$ that come
from writing in the above basis the generators of
$x_n^i\frb^{r+p^e-i}$ are divisible by $x_n$ if $i\geq p^e$, hence
map to zero in $S$. On the other hand, the generators coming from
$x_n^i\frb^{r+p^e-i}$ for $i\leq p^e-1$ are the same as the ones
obtained from writing the generators of $\frb^{r+p^e-i}$ in the
corresponding basis of $k\llbracket x_1,\ldots,x_{n-1}\rrbracket$
over $k^{p^e}\llbracket x_1^{p^e},\ldots,x_{n-1}^{p^e}\rrbracket$.
Moreover, it is clear that it is enough to consider only the largest
such ideal, namely $\frb^{r+1}$. This shows that
$\left(\fra^{r+p^e}\right)^{[1/p^e]}\cdot S=
\left((\fra/(x_n))^{r+1}\right)^{[1/p^e]}$, as claimed.
\end{proof}

\begin{corollary}\label{cor3}
If $R$, $S$ and $\fra$ are as in Proposition \ref{lem6}, and if $c>0$
is a jumping exponent for $\fra\cdot S$, then $c+r$ is a jumping
exponent for $\fra$.
\end{corollary}

We can give now the proof of our main result.

\begin{proof}[Proof of Theorem~\ref{thm1}]
For (1), suppose that we have a sequence of $F$-jumping exponents
$\{\alpha_m\}_m$ for $\fra$ having a finite accumulation point
$\alpha$. By Corollary~\ref{cor1}, we have $\alpha_m<\alpha$ for
$m\gg 0$. After replacing this sequence by a subsequence, we may
assume that $\alpha_m<\alpha_{m+1}$ for every $m$.

Let us write $R=R_1\times\cdots\times R_s$, where all $R_i$ are
domains. We have $\fra=\fra_1\times\cdots\times\fra_s$, and for
every $m$ there is $j$ such that $\alpha_m$ is a jumping exponent
for $\fra_j$. After replacing our sequence by a subsequence, we may
replace $R$ by some $R_j$ and therefore assume that $R$ is a domain.

By hypothesis, we can write $R\simeq S^{-1}(k[x_1,\ldots,x_n]/I)$
for some ideal $I$ and some multiplicative system $S$. Note that
$[k\colon k^p]<\infty$ (see Example~\ref{rem0}). Let us write $\fra=
S^{-1}(\frb/I)$, for some ideal $\frb \subset k[x_1, \dots, x_n]$ containing $I$.
 Note that $S^{-1}I$ is a prime ideal, hence it has
pure codimension, say  $r$, in
$S^{-1}k[x_1,\ldots,x_n]$. It follows from Corollary~\ref{cor3} that
$r+\alpha$ is an accumulation point for the jumping exponents of
$S^{-1}\frb$. Moreover, Proposition~\ref{prop4}(1) implies that
$r+\alpha$ is an accumulation point for the jumping exponents of
$\frb$. Therefore, in order to achieve a contradiction we may assume
that $R=k[x_1,\ldots,x_n]$.

Suppose now that $\fra$ is generated by polynomials of degree at
most $d$. In this case,  Proposition~\ref{prop6} implies that every
$\tau(\fra^{\alpha_m})$ is generated by polynomials of degree at
most $\lfloor\alpha d\rfloor$. Since the ideals
$\tau(\fra^{\alpha_m})$ form a strictly decreasing sequence of
ideals,  their  $k$-subvector spaces of polynomials of degree at most  $\lfloor \alpha d\rfloor$ form a strictly decreasing sequence of  subspaces of the finite dimensional vector space $k[x_1,\ldots,x_n]_{\leq\lfloor \alpha d\rfloor}.$ This contradiction
 completes the proof of (1).

Suppose now that $\alpha>0$ is an $F$-jumping exponent for $\fra$.
Proposition~\ref{lem5}(1) implies that all $p^e\alpha$ are
$F$-jumping exponents for $\fra$. We may assume that no $p^e\alpha$
is an integer, as otherwise $\alpha$ is clearly rational. Suppose
that $\fra$ is generated by $m$ elements and let $e_0$ be such that
$p^{e_0}\alpha>m$. We deduce from Proposition~\ref{lem5}(2) that
$\{p^e\alpha\}+m-1$ is a jumping exponent for every $e\geq e_0$.
Since all these numbers lie in $[m-1,m)$, it follows from (1) that
there are only finitely many such numbers. Therefore we can find
$e_1\neq e_2$ such that $p^{e_1}\alpha-p^{e_2}\alpha$ is an integer.
Hence $\alpha$ is a rational number.
\end{proof}

The ideas in the above proof can be used to explicitly estimate the
$F$-jumping exponents of an ideal in a polynomial ring in terms of
the degrees of its generators. Since we know that all $F$-jumping
exponents are rational, it is enough to bound their denominators.

\begin{proposition}\label{prop9}
Let $\fra\subseteq k[x_1,\ldots,x_n]$ be an ideal generated by $m$
polynomials of degree at most $d$. If $e_0$ is such that
$p^{e_0}>md$ and $N={{md+n}\choose{n}}$, then for every $F$-jumping
exponent $\alpha$ of $\fra$ we have $p^a(p^b-1)\alpha\in\NN$ for
some $a\leq e_0+N$ and $b\leq N$.
\end{proposition}

\begin{proof}
By Proposition~\ref{lem5}(2) it is enough to consider the case when
$\alpha\leq m$. Since $\tau(\fra^c)$ is generated by polynomials of
degree at most $md$ for every $c\leq m$, it follows that we have at
most $N=\dim_kk[x_1,\ldots,x_n]_{\leq md}$ jumping exponents of
$\fra$ in $(0,m]$.

We may assume that for every $e\leq e_0+N$ we have
$p^e\alpha\not\in\NN$, since otherwise our assertion is clear. In
particular, $\alpha>0$ and therefore  $\alpha\geq  1/d$ (note that
$\tau(\fra^c)=R$ for $c<1/d$ by Proposition~\ref{prop6}). Consider
now $e$, with $e_0\leq e\leq e_0+N$. Since $p^e\alpha\geq
\frac{p^{e_0}}{d}>m$, Proposition~\ref{lem5} implies that
$\{p^e\alpha\}+m-1$ is an $F$-jumping exponent of $\fra$. This gives
$N+1$ numbers in $(m-1,m)$ that are $F$-jumping exponents. We deduce
that there are $e_1<e_2$ in $\{e_0,\ldots,e_0+N\}$ such that
$p^{e_1}(p^{e_2-e_1}-1)\alpha\in\NN$, which completes the proof.
\end{proof}

\begin{remark}
We can get a bound independent of $m$ in the above proposition by
taking $m=n$. Indeed, note first that we may assume that $k$ is
infinite: otherwise take an infinite perfect extension $K$ of $k$,
and replace $\fra$ by $\fra\cdot K[x_1,\ldots,x_n]$ using
Remark~\ref{rem10}.

Now, for any maximal  ideal $\frm \subset k[x_1, \ldots, x_n] = R$,
it is clear (by Proposition \ref{prop4}) that the set of jumping
exponents for $\fra R_{\frm}$ is a subset of those for $\fra$. On
the other hand, for each of the finitely many jumping exponents
$\alpha$ in the range $[0, m]$, there is a  maximal ideal
$\frm_{\alpha}$ such that the inclusion $\tau(\fra^{c}) \subsetneqq
\tau(\fra^{c-\epsilon})$ is preserved under localization at
$\frm_{\alpha}$. It follows from Proposition \ref{lem5} (2)  that
the set of jumping exponents for $\fra$ is contained in (hence equal
to) the union of the sets of jumping exponents for the
$\frm_{\alpha}$, as we range over these finitely many
$\frm_{\alpha}$.

Now, fix generators $h_1,\ldots,h_m$ for $\fra$ of degree at most
$d$, and let $\frb $ be the ideal generated by
$g_i=\sum_{j=1}^ma_{i,j}h_j$ for $1\leq i\leq n$, where $a_{i,j}\in
k$.  It is well-known that for every maximal ideal $\frm$ in $R$,
the ideals $\frb R_{\frm} $ and $\fra R_{\frm}$ have the same
integral closure if the $a_{i,j}$ are general in $k$ (see, for
example \cite{positivity}, Example~9.6.19; note that in the proof
therein one does not use the assumption that the base field is
${\mathbb C}$).  Hence,  by choosing the $a_{i,j}$ sufficiently
general, we may assume that $\frb R_{\frm_{\alpha}} $ and $\fra
R_{\frm_{\alpha}}$ have the same integral closure, and hence the
same test ideals,  for each of our finitely many relevant maximal
ideals $\frm_{\alpha}$. It follows that the set of  jumping
exponents for $\fra$ is contained in the set of jumping exponents
for the $n$-generated ideal $\frb$. Thus we get the desired
statement for $\fra$ by applying Proposition~\ref{prop9} to $\frb$.
\end{remark}

\begin{remark}
Consider several ideals $\fra_1,\ldots,\fra_r$ in
$k[x_1,\ldots,x_n]$, where $k$ is $F$-finite. The same argument used
in the proof of Proposition~\ref{prop6} shows that if each $\fra_i$
is generated in degree at most $d_i$, then for every
$c_1,\ldots,c_r\in\RR_+$, the mixed test ideal
$\tau(\fra^{c_1}\cdots\fra^{c_r})$ can be generated in degree at
most $\sum_{i=1}^rc_id_i$. However, this assertion does not seem to
have such strong consequences in the case of several ideals. The
most optimistic expectation in this case is that for every
$b_1,\ldots,b_r$, the region
$$\{c\in\RR_+^n\mid c_i\leq b_i\,{\rm for}\,{\rm all}\,i\}$$
can be decomposed in a finite set of rational polytopes with
non-overlapping interiors, such that on the interior of each face of
such a polytope the mixed text ideal
$\tau(\fra_1^{c_1}\cdots\fra_r^{c_r})$ is constant.

If $k$ is contained in the algebraic closure of a finite field, then
we can find a finite field $k_0$ such that every $\fra_i$ is
generated by polynomials in $k_0[x_1,\ldots,x_n]$. In this case, it
follows that every $\tau(\fra_1^{c_1}\cdots\fra_r^{c_r})$ is
generated by polynomials with coefficients in $k_0$. In particular,
if the $c_i$ are bounded above, then the generators of the above
test ideal lie in a finite-dimensional vector space over the finite
field $k_0$. Therefore given $\fra_1,\ldots,\fra_r$, we can have
only finitely many possible test ideals
$\tau(\fra^{c_1}\cdots\fra^{c_r})$ if we bound the $c_i$.
\end{remark}


\end{document}